\documentstyle[11pt,twoside]{article}

\input{mssymb}
\makeatother

\newcommand{\mysec}[2]{%
\section*{\normalsize\hfil\sc {#1}. {#2}\hfill}%
\setcounter{theo}{0}\setcounter{equation}{0}\setcounter{section}{#1}%
\typeout{#1. #2}
\noindent}

\newcommand{\Proof}{\par\noindent{\em Proof. }}
\newcommand{\eop}{\nopagebreak\hspace*{\fill}$\Box$}

\newtheorem{theo}{Theorem}[section]
\newtheorem{lemma}[theo]{Lemma}

\newtheorem{prop}[theo]{Proposition}
\newtheorem{definition}[theo]{Definition}

\newcounter{abc}   
\newcounter{iiiii} 

\newenvironment{aequivalenz}
{\setcounter{iiiii}{0}
\begin{list}%
{{\rm (\roman{iiiii})}}
{\usecounter{iiiii}
\parsep=0pt plus 1pt
\topsep=1pt plus 2pt minus 1pt
\itemsep=1pt plus 2pt minus 1pt
\leftmargin=3\baselineskip
\labelsep=.6\baselineskip
\labelwidth=2.4\baselineskip
\rightmargin 0pt}%
}
{\end{list}}

\newenvironment{statements}%
{\setcounter{abc}{0}
\begin{list}%
{{\rm (\alph{abc})}}
{\usecounter{abc}
\parsep=0pt plus 1pt
\topsep=1pt plus 2pt minus 1pt
\itemsep=1pt plus 2pt minus 1pt
\leftmargin=3\baselineskip
\labelsep=.6\baselineskip
\labelwidth=2.4\baselineskip
\rightmargin 0pt}%
}
{\end{list}}

\newenvironment{punkte}%
{\begin{list}%
{$\bullet$}
{\leftmargin=3\baselineskip
\labelsep=\baselineskip
\labelwidth=2.5\baselineskip
\parsep=0pt plus 1pt
\topsep=1pt plus 2pt minus 1pt
\itemsep=1pt plus 2pt minus 1pt
\rightmargin 0pt}%
}
{\end{list}}

\newif\ifrefsc
\let\thebibliographyalt=\thebibliography                                %
\def\thebibliography#1                                                  %
 {\ifrefsc
 \section*{\normalsize\hfil\sc References\hfill}
 \small
 \list                    %
 {[\arabic{enumi}]}{\settowidth\labelwidth{[#1]}\leftmargin\labelwidth  %
 \advance\leftmargin\labelsep                                           %
 \usecounter{enumi}}                                                    %
 \def\newblock{\hskip .11em plus .33em minus .07em}                     %
 \sloppy\clubpenalty4000\widowpenalty4000\vfuzz 2pt                      %
 \sfcode`\.=1000\relax                                                  %
 \else\thebibliographyalt{#1}\fi}                                         %

\makeatletter
\def\ps@smallheadings{            
  \def\@oddfoot{}                 
  \def\@evenfoot{}                
  \def\@evenhead{\hbox to \textwidth {%
  \vbox{\hbox to \textwidth 
        {\thepage\hfil\strut{\footnotesize\sc\Autor}\hfil}\vss}}}
  \def\@oddhead{\hbox to \textwidth {%
\vbox{\hbox to \textwidth 
        {\hfil{\footnotesize\sc\Kurztitel}\hfil\strut \thepage}\vss}}}
    }
\makeatother

\makeatletter

\def\eqalignno#1{\displ@y \tabskip\@centering
  \halign to\displaywidth{\hfil$\@lign\displaystyle{##}~$\tabskip\z@skip
    &$\@lign\displaystyle{{}##}$\hfil\tabskip\@centering
    &\llap{$\@lign##$}\tabskip\z@skip\crcr
    #1\crcr}}

\makeatother

\def\ersteSeite{\vspace*{32pt plus 2pt minus 2pt}\begin{center}
{\Large\sf\Titel}\\[15pt]{\sc\Autor}\\[26pt plus 2pt minus 2pt]
\end{center}\Abstrakt\vspace{0pt plus 2pt }\thispagestyle{empty}}

\def\Abstrakt{\begin{quote}\small\noindent{\sc Abstract.}
\Abstrakttext\end{quote}}

\newcommand{\N}{{\Bbb N}}

\newcommand{\KK}{{\Bbb K}}

\newcommand{\Id}{I\mkern-1mud}

\newcommand{\alp}{\alpha}

\newcommand{\del}{\delta}

\newcommand{\lam}{\lambda}
\newcommand{\Lam}{\Lambda}
\newcommand{\sig}{\sigma}
\newcommand{\Sig}{\Sigma}

\newcommand{\klam}{\bigl(}
\newcommand{\mer}{\bigr)}

\newcommand{\loglike}[1]{\mathop{\rm #1}\nolimits}

\newcommand{\ex}{\loglike{ex}}

\newcommand{\lin}{\loglike{lin}}

\newcommand{\ran}{\loglike{ran}}

\newcommand{\supp}{\loglike{supp}}

\def\iinf{\mathop{\rm \vphantom{p}\inf}} 
 
\def\limp{\mathop{\rm \vphantom{p}\lim}} 

\newcommand{\dst}{\displaystyle}

\newcommand{\qqfa}{\qquad\forall}

\newcommand{\bea}{\begin{eqnarray*}}
\newcommand{\eea}{\end{eqnarray*}}
\newcommand{\beq}{\begin{equation}}
\newcommand{\eeq}{\end{equation}}
\newcommand{\begsta}{\begin{statements}}
\def\endsta{\end{statements}}
\newcommand{\begaeq}{\begin{aequivalenz}}
\def\endaeq{\end{aequivalenz}}

\newcommand{\iy}{\infty}
\newcommand{\dopu}{{:}\allowbreak\ }

\newcommand{\pel}{Pe{\l}czy\'{n}ski}

\newcommand{\rest}[2]{#1\raisebox{-0.3ex}{\mbox{$\mid_{#2}$}}} 

\def\Einzelzeile{(2.2)}

\refsctrue

\begin{document}


\def\Titel{ The Daugavet equation for operators\\[3pt]
 not fixing a copy of $C(S)$}
\def\Kurztitel{ Operators not fixing a copy of $C(S)$}
\def\Autor{Lutz Weis and Dirk Werner}
\def\Abstrakttext{
We prove the norm identity $\|\Id+T\| = 1+\|T\|$, which is known as
the Daugavet equation, for operators $T$ on $C(S)$ not fixing a copy
of $C(S)$, where $S$ is a compact metric space without isolated
points.
}

\ersteSeite

\mysec{1}{Introduction}%
An operator $T\dopu X\to X$ on a Banach space is said to satisfy the
{\em Daugavet equation\/} if
\beq\label{eq0}
\|\Id + T\| = 1 + \|T\|;
\eeq
this terminology is derived from Daugavet's theorem that a compact
operator on $C[0,1]$ satisfies (\ref{eq0}). Many authors have
established the Daugavet equation for various classes of operators,
e.g., the weakly compact ones, on various spaces; we refer to
\cite{Abra2}, \cite{AbraAB}, \cite{Kadets}, \cite{PliPop},
\cite{Dirk8}, \cite{Dirk10}, \cite{Woj92} and the references in these
papers for more information.

The most far-reaching result for operators on $L^{1}[0,1]$ is
due to Plichko and Popov \cite[Th.~9.2 and Th.~9.8]{PliPop} who
prove that an operator on $L^{1}[0,1]$ which does not fix a copy of
$L^{1}[0,1]$ satisfies the Daugavet equation. (As usual, $T\dopu X\to
X$ fixes a copy of a Banach space~$E$  if there is a subspace
$F\subset X$  isomorphic to $E$ such that $\rest{T}{F}$ is an (into-)
isomorphism.) 
In this paper we shall establish the corresponding result for
operators on spaces of continuous functions.

\begin{theo}\label{T1}
Let $(S,d)$ be a compact metric space without isolated points. If $T\dopu
C(S)\to C(S)$ does not fix a copy of $C(S)$, then $T$ satisfies the Daugavet
equation.
\end{theo}

In the next section we shall provide a direct proof of this theorem
which is measure theoretic in spirit, and in Section~3 we shall
give another argument which relies on a deep result due to
Rosenthal \cite{Ros-nons}. Though this proof is much shorter it is
less revealing than the one in Section~2; therefore we found it
worthwhile to present both arguments.

Either of our approaches depends on the analysis of the representing
kernel of an operator $T\dopu C(S)\to C(S)$; this is the family of
Borel measures on $S$ defined by $\mu_{s}= T^{*}\del_{s}$, $s\in S$.
In \cite{Dirk8} (see also \cite{Dirk10}) a necessary and sufficient
condition on the kernel of $T$ was established in order that $T$
satisfy the Daugavet equation. Here we record a special case which 
suffices for our needs.

\begin{lemma}\label{L1}
If the representing kernel $(\mu_{s})_{s\in S}$
 of a bounded linear operator $T\dopu
C(S)\to C(S)$  satisfies
\beq\label{eq1}
\inf_{s\in U} |\mu_{s}(\{s\})| =0
\eeq
for all nonvoid open subsets $U$ of $S$, then $T$ satisfies the Daugavet
equation.
\end{lemma}

Our results are valid for real and complex spaces; we denote the
scalar field by~$\KK$.

\mysec{2}{Proof of Theorem 1.1}%
We fix a diffuse probability measure $\lam$ on $S$ whose support
is~$S$. To show the existence of such a measure, consider a countable
basis of open sets $O_{n}\subset S$. Since each $O_{n}$ is perfect,
we know from a theorem of Bessaga and \pel\ (see \cite[p.~52]{Lac})
that there are
diffuse probability measures~$\lam_{n}$ with $\supp \lam_{n} \subset
O_{n}$; it remains to define $\lam = \sum 2^{-n} \lam_{n}$.

Now let $T\dopu C(S) \to C(S)$ be a bounded linear operator and
$(\mu_{s})_{s\in S}$ its representing kernel. It is enough to
show that $T$ fixes a copy of $C(S)$ if (\ref{eq1}) of
Lemma~\ref{L1} fails for some open set~$U$. Thus we assume, for some
open $U\subset S$ and $\alp>0$, that
\beq\label{eq2}
|\mu_{s}(\{s\})| \ge \alp \qquad\forall s\in U.
\eeq
We shall also assume without loss of generality that $\|T\|=1$.

We shall decompose the $\mu_{s}$ into their atomic and diffuse parts.
It is instrumental for our argument that this can be accomplished in
a measurable fashion, as was proved by Kalton \cite[Th.~2.10]{Kal-Endo}. 
More precisely we can write
$$
\mu_{s}= \sum_{n=1}^{\iy} a_n(s) \del_{\sig_{n}(s)} +  \nu_{s}
$$
where
\smallskip
\begin{punkte}
\item
each $a_{n}\dopu S \to \KK$ is measurable for the
completion $\Sig_{\lam}$ of the Borel sets of $S$ with respect to $\lam$,
\item
each $\sig_{n}\dopu S \to S$ is $\Sig_{\lam}$-Borel-measurable,
\item
each $\nu_{s}$ is diffuse, and $s\mapsto |\nu_{s}|$
is $\Sig_{\lam}$-Borel-measurable (we are considering the weak$^{*}$ Borel
sets of the unit ball of $C(S)^{*}$),
\item
$|a_{1}(s)| \ge |a_{2}(s)| \ge \ldots$ for all $s\in S$,
\item
$\sig_{n}(s)\neq \sig_m(s)$ for all $s\in S$ whenever $n\neq m$,
\item
$\sum_{n=1}^{\iy} |a_{n}(s)| \le1$ for all $s\in S$.
\end{punkte}
\smallskip

Put $\beta = \lam(U)>0$. Applying the Egorov and Lusin theorems, we
may find a compact subset $S_{1}\subset S$ with $\lam(S_{1}) \ge
1-\beta/2$ such that
\smallskip
\begin{punkte}
\item
each $\rest{a_{n}}{S_{1}}$ is continuous on $S_{1}$,
\item
each $\rest{\sig_{n}}{S_{1}}$ is continuous on $S_{1}$,
\item
$\sum |a_{n}(s)|$ converges uniformly on $S_{1}$,
\item
$\dst \limp_{n\to\iy} \sup_{s\in S_{1}} |\nu_{s}|\klam B(t,1/n) \mer \to 0$
for all $t\in S$; \hfill\Einzelzeile
\addtocounter{equation}{1}%
\end{punkte}
\smallskip\par\noindent
here $B(t,r)$ denotes the closed ball with centre $t$ and radius~$r$.
To prove \Einzelzeile\  
fix a countable dense subset $Q$ of $S$. Let ${\bf U}_{n}$
be the countable set of all finite covers of $S$     by means of open
balls of radius~$1/n$ with centres in~$Q$. Then
$$
\varphi_{n}(s) = \iinf_{{\cal C}\in {\bf U}_{n}} \sup_{B\in{\cal C}}
|\nu_{s}|(B)
$$
defines a sequence of $\Sig_{\lam}$-measurable functions, and we have
$\varphi_{n}\to 0$ pointwise, since the $\nu_{s}$ are diffuse and $S$
is compact. By Egorov's theorem $(\varphi_{n})$ tends to $0$
uniformly on a large subset of $S$. This easily implies our claim.

After discarding relatively isolated points of
$S_{1}$, if necessary,  
we come up with a perfect subset $S_{1}$ with the above
properties; note that $U\cap S_{1} \neq \emptyset$ since this
intersection has positive $\lam$-measure. Next we pick $N\in\N$ such
that
\beq\label{eq4}
\sum_{n>N} |a_{n}(s)| \le \alp/3 \qquad \forall s\in S_{1}.
\eeq

Let $s\in U\cap S_{1}$. Then
$$
\alp \le |\mu_{s}(\{s\})| = \left| \sum_{n=1}^{\iy} a_{n}(s)
\del_{\sig_{n}(s)}(\{s\}) \right|.
$$
Consequently $\sig_{k}(s) = s$ for some $k$; in addition, we must
have 
\beq\label{eq5}
|a_{k}(s)|\ge\alp
\eeq 
for this~$k$ so that $k\le N$. This shows
that
$$
U\cap S_{1} = \bigcup_{l=1}^{N} \{s\in U\cap S_{1}\dopu
\sig_{l}(s)=s\}.
$$
Each of these $N$ sets is relatively closed in $U\cap S_{1}$, so one
of them, with index $k$ say, contains a relatively interior
point~$s_{0}$. Therefore, there exists some $\del>0$ with the
following properties:

For $B_{l} = \sig_{l}\klam S_{1}\cap B(s_{0},\del) \mer$, $l=1,\ldots,
N$, we have
\smallskip
\begsta
\item
$B_{k}= S_{1}\cap B(s_{0},\del)\subset S_{1}\cap U$, 
and $\sig_{k}(s)=s$ for all $s\in
B_{k}$,
\item
$B_{l}\cap B(s_{0},\del) =\emptyset$ for $l\neq k$ 
(since $s_{0}= \sig_{k}(s_{0})
\neq \sig_{l}(s_{0})$ and these functions are continuous on $S_{1}$),
\item
$|\nu_{s}| \klam  B(s_{0},\del) \mer \le \alp/3$ for all $s\in S_{1}$
(see~\Einzelzeile).
\endsta
\smallskip

Now let $V= S_{1} \cap B(s_{0},\del/2)$.  This is a compact perfect
subset of $S$ and thus $C(V)$ is isomorphic to $C(S)$  by Milutin's
theorem \cite[Th.~III.D.19]{Woj-Buch}. Also let
$W= \{s\in S\dopu d(s,s_{0})\ge \del\}$.
Again, this is a compact subset of $S$, $V\cap W=\emptyset$, and we
have from (a) and (b) above
\smallskip
\begsta
\item[\rm (a$'$)]
$\sig_{k}(s)=s$ for all $s\in V$,
\item[\rm (b$'$)]
$B_{l}\subset W$ for all $l=1,\ldots,N$, $l\neq k$.
\endsta
\smallskip\par\noindent
By the Borsuk-Dugundji theorem (see e.g.\
\cite[Cor.~III.D.17]{Woj-Buch}) there exists an isometric linear
extension operator $L\dopu C(V\cup W) \to C(S)$. Let $E=\{ f\in
C(V\cup W)\dopu \rest{f}{W}=0\}$ and $F=L(E)$; thus $F\cong E \cong
C(V)$, and $F$ is isomorphic to~$C(S)$.

We finally prove that $\rest{T}{F}$ is an isomorphism. In fact, let
us show that
$$
\|Tf\| \ge \frac\alp3 \, \|f\| \qquad\forall f\in F.
$$
Suppose $f\in F$ with $\|f\|=1$. Then $\|\rest{f}{V}\|=1$ and
$\rest{f}{W}=0$. Pick $s\in V$ with $|f(s)|=1$. Then we have (recall
that $V\subset S_{1} \cap U$)
\bea
\|Tf\| &\ge& |Tf(s)| \\
&=&
\left| \sum_{n=1}^{\iy} a_n(s) f\klam \sig_{n}(s) \mer +
   \int_{S} f\, d\nu_{s} \right|  \\
&\ge&
|a_{k}(s)| - \sum_{n>N} |a_{n}(s)| - \int_{S} |f|\,d|\nu_{s}| \\
&& \qquad\mbox{(by (a$'$) and (b$'$))}  \\
&\ge&
\alp - \frac\alp3 - |\nu_{s}|\klam B(s_{0},\del) \mer \\
&& \qquad\mbox{(by (\ref{eq5}), (\ref{eq4}) and since $\rest{f}{W}=0$)} \\
&\ge& \frac\alp3 \qquad\mbox{(by (c))}.
\eea
This completes the proof of the theorem.
\eop

\mysec{3}{Operators whose adjoints have separable ranges}%
In this section we discuss another, formally shorter, approach to
Theorem~\ref{T1} via the following proposition.

\begin{prop}\label{P3.1}
Suppose $S$ is a compact Hausdorff space without isolated points. If
$T\dopu C(S)\to C(S)$  is a bounded linear operator such that
$\ran(T^{*})$ is separable, then $T$ satisfies the Daugavet equation.
\end{prop}

\Proof
Let $(\mu_{s})_{s\in S}$ be the representing kernel of $T$, and let
$\{\mu_{s_{n}}\dopu n\in\N\}$ be dense in
$\{\mu_{s}\dopu s\in S\}$. We define $\mu=\sum 2^{-n} |\mu_{s_{n}}|$;
then all the $\mu_{s}$ are absolutely continuous with respect
to~$\mu$. Hence $\mu_{s}(\{s\})\neq 0 $ only if $|\mu|(\{s\})\neq0$.
Therefore, $\{s\in S\dopu \mu_{s}(\{s\}) \neq0\}$ is countable, and,
since $S$ is perfect, its complement is dense by Baire's theorem.
This implies that the condition of Lemma~\ref{L1} is satisfied, and
the proposition is proved.
\eop
\bigskip

Proposition~\ref{P3.1} obviously covers the case of compact
operators, operators factoring through a space with a separable dual, and, in
the case of metric~$S$, weakly compact operators; for the latter
observe that $T^{**}\klam C(S)^{**} \mer \subset C(S)$ so that $T^{**}$
and, consequently, $T^{*}$ have separable ranges. But even for
nonmetrizable perfect~$S$ the Daugavet equation for weakly compact
operators on $C(S)$ can be derived from the argument of
Proposition~\ref{P3.1}; one only has to recall that by a theorem of
Bartle, Dunford and Schwartz \cite[p.~306]{DS} 
also in this case all the $\mu_{s}$
are absolutely continuous with respect to some finite measure~$\mu$.
This seems to yield the easiest proof of the Daugavet equation for
these classes of operators. 

If $S$ is a metric space then Proposition~\ref{P3.1} clearly follows
from Theorem~\ref{T1}. On the other hand, 
Rosenthal \cite{Ros-nons} (see also \cite{Weis4}) 
has shown that an operator $T\dopu C(S)\to
C(S)$, $S$ an uncountable 
compact metric space, with $\ran(T^{*})$ nonseparable fixes a
copy of $C(S)$. Therefore Proposition~\ref{P3.1} implies
Theorem~\ref{T1} so that both these results are in fact equivalent.
Moreover we note that in this setting $T$ does not fix a
copy of $C(S)$ if and only if $T$ does not fix a copy of $\ell^{1}$
if and only if $T$ is weakly sequentially precompact,
i.e., if $(f_{n})$ is a bounded sequence in $C(S)$, then $(Tf_{n})$
admits a weak Cauchy subsequence. 

We finally extend Proposition~\ref{P3.1} to the class of nicely
embedded Banach spaces introduced in \cite{Dirk10} 
so that we obtain a unified
approach to several results of that paper.

Let $S$ be a Hausdorff
topological space, and let
$C^{b}(S)$ be the sup-normed Banach space of all bounded continuous
scalar-valued functions. The functional $f\mapsto f(s)$ on $C^{b}(S)$
is denoted by $\del_{s}$. We say that a linear map $J\dopu X \to
C^{b}(S)$ on a Banach space $X$ is a {\em nice embedding\/} and that
$X$ is {\em nicely embedded\/} into $C^{b}(S)$ if $J$ is an isometry
such that for all $s\in S$ the following properties hold:
\smallskip
\begsta\em
\item[\rm(N1)]
For $p_{s}:= J^{*}(\del_{s})\in X^{*}$ we have $\|p_{s}\|=1$.
\item[\rm(N2)]
$\lin\{p_{s}\}$ is an $L$-summand in $X^{*}$.
\endsta
\smallskip
The latter condition means that there are
projections $\Pi_{s}$ from $X^{*}$ onto
$\lin\{p_{s}\}$ such that
$$
\|x^{*}\| = \|\Pi_{s}(x^{*})\| + \|x^{*}-\Pi_{s}(x^{*})\| \qqfa x^{*}\in
X^{*}.
$$

We will also need the equivalence relation
$$
s \sim t \ \mbox{ if and only if }\  \Pi_{s} =\Pi_{t}
$$
on $S$. Then $s$ and $t$ are equivalent if and only if $p_{s}$ and
$p_{t}$  are linearly dependent, which implies by (N1) that
$p_{t}=\lam p_{s}$ for some scalar of modulus~1. The equivalence
classes of this relation are obviously closed.

We will consider the following nondiscreteness condition.
\smallskip
\begsta\em
\item[\rm (N3)]
None of the equivalence classes $Q_{s}=\{t\in S\dopu s\sim
t\}$ contains an interior point.
\endsta
\smallskip
If the set $\{p_{s}\dopu s\in S\}$ is linearly independent, this
simply means:
\smallskip
\begsta\em
\item[\rm (N$3'$)]
$S$ does not contain an isolated point.
\endsta
\smallskip
By (N2), the $p_{s}$ are linearly independent as soon as they are
pairwise linearly independent.

\begin{prop}\label{P3.2}
Let $S$ be a Baire topological Hausdorff space and suppose that $X$
is nicely embedded into $C^{b}(S)$ so that additionally\/ {\rm (N3)}
holds. If $T\dopu X\to X$ is a bounded linear operator with
$\ran(T^{*})$ separable, then $T$ satisfies the Daugavet equation.
\end{prop}

\Proof
We stick to the above notation and put $q_{s}= T^{*}(p_{s})$. 
By \cite[Prop.~2.1]{Dirk10} it is sufficient
to prove that
$$
S' = \{t\in S\dopu \ \Pi_{t}(q_{s}) =0 \ \forall s\}
$$
is dense in $S$. If $\{q_{s_{n}}\dopu n\in \N\}$ denotes a dense
subset of $\{q_{s}\dopu s\in S\}$, then $t\notin S'$ if and only if
$\Pi_{t}(q_{s_{n}})\neq0$ for some~$n$.  This implies that
$S\setminus S'$ consists of at most countably many equivalence
classes for the equivalence relation~${\sim}$ (cf.\ 
\cite[Lemma~2.3]{Dirk10}),
and by (N3) and Baire's theorem $S'$ must be dense.
\eop
\bigskip

It is proved in \cite{Dirk10} that the following classes of Banach spaces
satisfy the assumptions of Proposition~\ref{P3.2}:
\smallskip
\begin{punkte}
\item
$X$ is a function algebra whose Choquet boundary does not contain
isolated points,
\item
$X$ is a real $L^{1}$-predual space so that the set $\ex B_{X^{*}}$
of extreme functionals does not contain isolated points (for the
weak$^{*}$ topology),
\item
$X$ is a complex $L^{1}$-predual space for which the quotient space
 $\ex B_{X^{*}}/{\sim}$ does not contain isolated points, where $\sim$
means linear dependence,
\item
$X$ is a space of type $C_{\Lam}$ for certain subsets $\Lam$ of
abelian discrete groups.
\end{punkte}

%
%
%
%
\typeout{References}

%
%
%
\small
\bigskip
\noindent
Mathematisches Institut~I, Universit\"at Karlsruhe,\\
D-76\,128 Karlsruhe, Germany; \ e-mail:
Lutz.Weis@math.uni-karlsruhe.de

\smallskip\noindent
I.~Mathematisches Institut, Freie Universit\"at Berlin,
Arnimallee 2--6, \\
D-14\,195 Berlin, Germany; \ 
e-mail: 
werner@math.fu-berlin.de

\end{document}